\documentclass{article}
\usepackage{amsfonts}

\newtheorem{theorem}{Theorem}[section]

\newtheorem{proposition}[theorem]{Proposition}

\newenvironment{proof}[1][Proof]{\noindent\textbf{#1.} }{\ \rule{0.5em}{0.5em}}

\begin{document}

\title{Differential equations driven by H\"{o}lder continuous functions of
order greater than $1/2$}
\author{Yaozhong Hu
\thanks{%
Y. Hu is supported in part by the National Science Foundation under Grant
No. DMS0204613 and DMS0504783 }  \ and \  David Nualart   \\
Department of Mathematics\thinspace ,\ University of Kansas\\
405 Snow Hall\thinspace ,\ Lawrence, Kansas 66045-2142\\
}
\date{}
\maketitle

\begin{abstract}
We derive estimates for the solutions to differential equations driven
by a H\"older continuous function of order $\beta>1/2$. 
As an application we deduce the existence of moments for the solutions to
stochastic partial differential equations driven by a fractional Brownian
motion with Hurst parameter $H>\frac{1}{2}$. 
\end{abstract}

\section{Introduction}

We are interested in the solutions of differential equations on $\mathbb{R}%
^{m}$ of the form 
\begin{equation}
x_{t}=x_{0}+\int_{0}^{t}f(x_{r})dy_{r},  \label{0.1}
\end{equation}%
where the driving force  $y:\mathbb{[}0,\infty )\rightarrow \mathbb{R}^{m}$
is a H\"{o}lder continuous function of order ${\beta }>1/2$. \ If the
function  $f:\mathbb{R}^{d}\rightarrow \mathbb{R}^{md}$ has bounded partial
derivatives which are H\"{o}lder continuous of order $\lambda >\frac{1}{%
\beta }-1$, then there is a unique  solution $x:\mathbb{R}^{m}\rightarrow 
\mathbb{R}$, which has bounded $\frac{1}{\beta }$-variation on any finite
interval.   These results have been proved by Lyons in \cite{Ly} using the $p$%
-variation norm and the technique introduced by Young in \cite{Yo}.  The
integral appearing in (\ref{0.1}) is then a Riemann-Stieltjes integral. 

In \cite{Za} Z\"ahle has introduced a generalized Stieltjes integral using
the techniques of fractional calculus. This integral is expressed in terms
of  fractional derivative operators and it coincides with the
Riemann-Stieltjes integral $\int_{0}^{T}fdg$, when the functions $f$ and $g$
are  H\"{o}lder continuous of orders $\lambda $ and $\mu $, respectively and 
$\lambda +\mu >1$ (see Proposition \ref{p1} below). Using this formula for
the Riemann-Stieltjes integral, Nualart and  R\u{a}\c{s}canu have    
obtained  in  \cite{NR} the existence of a unique  solution for a class
of general differential equations that includes (\ref{0.1}). Also they have
proved that the solution of (\ref{0.1}) is bounded on a finite interval $[0,T]$ by $%
C_{1}\exp (C_{2}\left\Vert y\right\Vert _{0,T,\beta }^{\kappa })$, where $%
\kappa >\frac{1}{\beta }$ if $f$ is bounded and $\kappa >\frac{1}{1-2\beta }$
is $f$ has linear growth. Here $\Vert y\Vert _{{0,T,\beta }}$ denotes the $\beta $%
-H\"{o}lder norm of $y$ on the time interval $[0,T]$. These estimates are
based on a suitable application of Gronwall's lemma. It turns out that the
estimate in the linear growth case is unsatisfactory because $\kappa $ tends
to infinity as $\beta $ tends to $1/2$. \ 

The main purpose of this paper is to obtain better estimates for the solution $x_t$ 
in the case
where $f$ is bounded or   has linear growth using a direct approach based
on formula (\ref{1.8}). In the case where $f$ is bounded we  estimate $%
\sup_{0\leq t\leq T}|x_{t}|$ by%
\[
C\left( 1+\Vert y\Vert _{0,T,{\beta }}^{\frac{1}{{\beta }}}\right) 
\]%
and if $f$ has linear growth we obtain the estimate%
\[
C_{1}\exp \left( C_{2}\Vert y\Vert _{0,T,{\beta }}^{\frac{1}{{\beta }}%
}\right) .
\]%
In Theorem \ref{t1} we provide explicit dependence on $f$ and $T$   for the constants $C$, $%
C_{1}$ and $C_{2}$.

Another novelty of this paper is that we establish    
 the explicit dependence of the solution $x_t$  to  (\ref{0.1})  on the initial condition $x_0$,
the driving control $y$ and the coefficient $f$ (Theorem 3.2).  Similar results
are obtained for the  case $1/3<\beta<1/2$ in a forthcoming  paper \cite{HN}. 

As an application we deduce the existence of moments for the solutions to
stochastic partial differential equations driven by a fractional Brownian
motion with Hurst parameter $H>\frac{1}{2}$.  We also discuss \ the
regularity of the solution in the sense of Malliavin Calculus, improving the
results of Nualart and Saussereau \cite{NS}, and we apply the techniques of
the Malliavin calculus to establish the existence of densities under
suitable   non-degeneracy conditions. 

 \setcounter{equation}{0}
\section{Fractional integrals and derivatives}

Let $a,b\in \mathbb{R}$ with $a<b$. Let $f\in L^{1}\left( a,b\right) $ and $%
\alpha >0.$ The left-sided and right-sided fractional Riemann-Liouville
integrals of $f$ of order $\alpha $ are defined for almost all $x\in \left(
a,b\right) $ by 
\[
I_{a+}^{\alpha }f\left( t\right) =\frac{1}{\Gamma \left( \alpha \right) }%
\int_{a}^{t}\left( t-s\right) ^{\alpha -1}f\left( s\right) ds 
\]%
and 
\[
I_{b-}^{\alpha }f\left( t\right) =\frac{\left( -1\right) ^{-\alpha }}{\Gamma
\left( \alpha \right) }\int_{t}^{b}\left( s-t\right) ^{\alpha -1}f\left(
s\right) ds, 
\]%
respectively, where $\left( -1\right) ^{-\alpha }=e^{-i\pi \alpha }$ and $%
\Gamma \left( \alpha \right) =\int_{0}^{\infty }r^{\alpha -1}e^{-r}dr$ is
the Euler gamma function. Let $I_{a+}^{\alpha }(L^{p})$ (resp. $%
I_{b-}^{\alpha }(L^{p})$) be the image of $L^{p}(a,b)$ by the operator $%
I_{a+}^{\alpha }$ (resp. $I_{b-}^{\alpha }$). If $f\in I_{a+}^{\alpha
}\left( L^{p}\right) \ $ (resp. $f\in I_{b-}^{\alpha }\left( L^{p}\right) $)
and $0<\alpha <1$ then the Weyl derivatives are defined as 
\begin{equation}
D_{a+}^{\alpha }f\left( t\right) =\frac{1}{\Gamma \left( 1-\alpha \right) }%
\left( \frac{f\left( t\right) }{\left( t-a\right) ^{\alpha }}+\alpha
\int_{a}^{t}\frac{f\left( t\right) -f\left( s\right) }{\left( t-s\right)
^{\alpha +1}}ds\right)  \label{1.1}
\end{equation}%
and 
\begin{equation}
D_{b-}^{\alpha }f\left( t\right) =\frac{\left( -1\right) ^{\alpha }}{\Gamma
\left( 1-\alpha \right) }\left( \frac{f\left( t\right) }{\left( b-t\right)
^{\alpha }}+\alpha \int_{t}^{b}\frac{f\left( t\right) -f\left( s\right) }{%
\left( s-t\right) ^{\alpha +1}}ds\right)  \label{1.2}
\end{equation}%
where $a\leq t\leq b$ (the convergence of the integrals at the singularity $%
s=t$ holds point-wise for almost all $t\in \left( a,b\right) $ if $p=1$ and
moreover in $L^{p}$-sense if $1<p<\infty $).

For any $\lambda \in (0,1)$, we denote by $C^{\lambda }(a,b)$ the space of $%
\lambda $-H\"{o}lder continuous functions on the interval $[a,b]$.   We will make use of the notation%
\[
\left\| x\right\| _{a,b,\beta }=\sup_{a\leq \theta <r\leq b}\frac{%
|x_{r}-x_{\theta }|}{|r-\theta |^{\beta }}, 
\]%
and%
\[
\Vert x||_{a,b,\infty }=\sup_{a\leq r\leq b}|x_{r}|, 
\]%
where $x:\mathbb{R}^{d}\rightarrow \mathbb{R}$ is a given continuous
function.

Recall
from \cite{SKM} that we have:

\begin{itemize}
\item If $\alpha <\frac{1}{p}$ and $q=\frac{p}{1-\alpha p}$ then 
\[
I_{a+}^{\alpha }\left( L^{p}\right) =I_{b-}^{\alpha }\left( L^{p}\right)
\subset L^{q}\left( a,b\right) . 
\]
\item If $\alpha >\frac{1}{p}$ then%
\[
I_{a+}^{\alpha }\left( L^{p}\right) \,\cup \,I_{b-}^{\alpha }\left(
L^{p}\right) \subset C^{\alpha -\frac{1}{p}}\left( a,b\right) . 
\]%
\end{itemize}
The following inversion formulas hold: 
\begin{eqnarray}
I_{a+}^{\alpha }\left( D_{a+}^{\alpha }f\right) &=&f,\quad \quad \;\forall
f\in I_{a+}^{\alpha }\left( L^{p}\right)  \label{1.4} \\
I_{a-}^{\alpha }\left( D_{a-}^{\alpha }f\right) &=&f,\quad \quad \;\forall
f\in I_{a-}^{\alpha }\left( L^{p}\right)  \label{1.3}
\end{eqnarray}%
and 
\begin{equation}
D_{a+}^{\alpha }\left( I_{a+}^{\alpha }f\right) =f,\quad D_{a-}^{\alpha
}\left( I_{a-}^{\alpha }f\right) =f,\quad \;\forall f\in L^{1}\left(
a,b\right) \,.  \label{1.5}
\end{equation}%
On the other hand, for any $f,g\in L^{1}(a,b)$ we have 
\begin{equation}
\int_{a}^{b}I_{a+}^{\alpha }f(t)g(t)dt=(-1)^{\alpha
}\int_{a}^{b}f(t)I_{b-}^{\alpha }g(t)dt\,,  \label{1.6}
\end{equation}%
and for $f\ \in I_{a+}^{\alpha }\left( L^{p}\right) $ and $g\in
I_{a-}^{\alpha }\left( L^{p}\right) $ we have%
\begin{equation}
\int_{a}^{b}D_{a+}^{\alpha }f(t)g(t)dt=(-1)^{-\alpha
}\int_{a}^{b}f(t)D_{b-}^{\alpha }g(t)dt.  \label{1.7}
\end{equation}%
Suppose that $f\in C^{\lambda }(a,b)$ and $g\in C^{\mu }(a,b)$ with $\lambda
+\mu >1$. Then, from the classical paper by Young \cite{Yo}, the
Riemann-Stieltjes integral $\int_{a}^{b}fdg$ exists. The following
proposition can be regarded as a fractional integration by parts formula,
and provides an explicit expression for the integral $\int_{a}^{b}fdg$ in
terms of fractional derivatives (see \cite{Za}).

\begin{proposition}
\label{p1} Suppose that $f\in C^{\lambda }(a,b)$ and $g\in C^{\mu }(a,b)$
with $\lambda +\mu >1$. Let ${\lambda }>\alpha $ and $\mu >1-\alpha $. Then
the Riemann Stieltjes integral $\int_{a}^{b}fdg$ exists and it can be
expressed as%
\begin{equation}
\int_{a}^{b}fdg=(-1)^{\alpha }\int_{a}^{b}D_{a+}^{\alpha }f\left( t\right)
D_{b-}^{1-\alpha }g_{b-}\left( t\right) dt,  \label{1.8}
\end{equation}%
where $g_{b-}\left( t\right) =g\left( t\right) -g\left( b\right) $.
\end{proposition}

 \setcounter{equation}{0}
\section{Estimates for the solutions of differential equations}

Suppose that $y:\mathbb{[}0,\infty )\rightarrow \mathbb{R}^{m}$ is a H\"{o}%
lder continuous function of order ${\beta }>1/2$. Fix an initial condition $%
x_{0}\in \mathbb{R}^{d}$ and consider the following differential equation 
\begin{equation}
x_{t}=x_{0}+\int_{0}^{t}f(x_{r})dy_{r},  \label{2.1}
\end{equation}%
where $f:\mathbb{R}^{d}\rightarrow \mathbb{R}^{md}$ is given function. Lyons
has proved in \cite{Ly} that Equation (\ref{2.1}) has a unique solution if $f$
is continuously differentiable and it has a derivative $f^{\prime }$ which
is \ bounded and locally$\ $H\"{o}lder continuous of order $\lambda >\frac{1}{%
\beta }-1$.

Our aim is to obtain estimates on $x_{t}$ which are better than those given
by Nualart and R\u{a}\c{s}canu in \cite{NR}.  

\begin{theorem}
\label{t1} Let $f$ be a continuously differentiable such that $f^{\prime }$
is \ bounded and locally$\ $H\"{o}lder continuous of order $\lambda >\frac{1%
}{\beta }-1$.
\begin{description}
\item[(i)] Assume that $f$ is also bounded. Then, there is a constant $k$,
which depends only on $\beta $, such that for all $T$, 
\begin{equation}
\sup_{0\leq t\leq T}|x_{t}|\leq |x_{0}|+kT\Vert f\Vert _{\infty }\Vert
f^{\prime }\Vert _{\infty }^{\frac{1-{\beta }}{{\beta }}}\Vert y\Vert _{0,T,{%
\beta }}^{\frac{1}{{\beta }}}.  \label{2.4}
\end{equation}
\item[(ii)] Assume that $f$ satisfies the linear growth condition 
\begin{equation}
|f(x)|\leq a_{0}+a_{1}|x|,  \label{2.6}
\end{equation}%
where $a_{0}\geq 0$ and $\ a_{1}\geq 0$. Then there is a
constant $k$ depending only on $\beta $, such that for all $T$, 
\begin{equation}
\sup_{0\leq t\leq T}|x_{t}|\leq 2^{kT\left[ \Vert f^{\prime }\Vert _{\infty
}\vee a_{0}\vee a_{1}\right] ^{1/\beta }\Vert y\Vert _{0,T,{\beta }}^{1/{%
\beta }}}(|x_{0}|+1)\,.  \label{2.5}
\end{equation}
\end{description}
\end{theorem}

\begin{proof}
Without loss of generality we assume that $d=m=1$. Assume first that $f$ is
bounded. Set $\ $ $\Vert y\Vert _{{\beta }}=\Vert y\Vert _{0,T,{\beta }}.$
Let $\alpha >1/2$ such that $\ \alpha >1-\beta $. First we use the
fractional integration by parts formula given in Proposition \ref{p1} to
obtain for all $s,t\in \lbrack 0,T]$, 
\[
|\int_{s}^{t}f(x_{r})dy_{r}|\leq \int_{s}^{t}|D_{s+}^{\alpha
}f(x_{r})\,D_{t-}^{1-\alpha }y_{t-}(r)|dr. 
\]%
From (\ref{1.2})  and  (\ref{1.1})    it is easy to see 
\begin{equation}
|D_{t-}^{1-\alpha }y_{t-}(r)|\leq k\Vert y\Vert _{r,t,{\beta }%
}|t-r| ^{\alpha +{\beta }-1}\leq k\Vert y\Vert _{{\beta }}|t-r| ^{\alpha +{%
\beta }-1}  \label{2.2}
\end{equation}%
and 
\begin{equation}
|D_{s+}^{\alpha }f(x_{r})|\leq k\left[ \left\| f\right\| _{\infty
}(r-s)^{-\alpha }+\ \left\| f^{\prime }\right\| _{\infty }\Vert x\Vert _{s,t,%
{\beta }}(r-s)^{\beta -\alpha }\right].  \label{2.3}
\end{equation}%
Therefore 
\begin{eqnarray*}
|\int_{s}^{t}f(x_{r})dy_{r}| &\leq &k\Vert y\Vert _{{\beta }}\int_{s}^{t}%
\left[ \left\| f\right\| _{\infty }(r-s)^{-\alpha }(t-r)^{\alpha +{\beta }%
-1}\right. \\
&&\left. +\left\| f^{\prime }\right\| _{\infty }\Vert x\Vert _{s,t,{\beta }%
}(r-s)^{{\beta }-\alpha }(t-r)^{\alpha +{\beta }-1}\right] dr \\
&\leq &k\Vert y\Vert _{{\beta }}\left[ \left\| f\right\| _{\infty }(t-s)^{{%
\beta }}+\left\| f^{\prime }\right\| _{\infty }\Vert x\Vert _{s,t,{\beta }%
}(t-s)^{2{\beta }}\right] \,.
\end{eqnarray*}%
Consequently, we have%
\[
\Vert x\Vert _{s,t,{\beta }}\leq k\Vert y\Vert _{{\beta }}\left[ \left\|
f\right\| _{\infty }+\left\| f^{\prime }\right\| _{\infty }\Vert x\Vert
_{s,t,{\beta }}(t-s)^{{\beta }}\right] \,. 
\]%
Hence, %
\[
\Vert x\Vert _{s,t,{\beta }}\leq k\Vert y\Vert _{{\beta }}\left\| f\right\|
_{\infty }\left( 1-k\left\| f^{\prime }\right\| _{\infty }\Vert y\Vert _{{%
\beta }}(t-s)^{{\beta }}\right) ^{-1}. 
\]%
Therefore, 
\begin{eqnarray*}
\Vert x||_{s,t,\infty } &\leq &|x_{s}|+\Vert x\Vert _{s,t,{\beta }%
}(t-s)^{\beta } \\
&\leq &|x_{s}|+k\Vert y\Vert _{{\beta }}\left\| f\right\| _{\infty }\left(
1-k\left\| f^{\prime }\right\| _{\infty }\Vert y\Vert _{{\beta }}(t-s)^{{%
\beta }}\right) ^{-1}(t-s)^{\beta }.
\end{eqnarray*}%
Let $A:=k\left\| f^{\prime }\right\| _{\infty }\Vert y\Vert _{{\beta }}$.
Divide the interval $[0,T]$ into $n=T/\Delta $ subintervals and apply the
above inequality on the interval $[0,\Delta ]$, $[\Delta ,2\Delta ]$ and so
on recursively to obtain 
\[
\sup_{0\leq t\leq T}|x_{t}|\leq |x_{0}|+kT\left\| f\right\| _{\infty }\Vert
y\Vert _{{\beta }}(1-A\Delta ^{\beta })^{-1}\Delta ^{{\beta }-1}. 
\]%
With the choice $\Delta =\left( \frac{1-{\beta }}{A}\right) ^{\frac{1}{{%
\beta }}}$ we get 
\begin{eqnarray*}
\sup_{0\leq t\leq T}|x_{t}| &\leq &|x_{0}|+kT\left\| f\right\| _{\infty
}\Vert y\Vert _{{\beta }}\frac{1}{{\beta }(1-{\beta })^{\frac{1-{\beta }}{{%
\beta }}}}\left( k\left\| f^{\prime }\right\| _{\infty }\Vert y\Vert _{{%
\beta }}\right) ^{\frac{1-{\beta }}{{\beta }}} \\
&=&|x_{0}|+kT\left\| f^{\prime }\right\| _{\infty }^{\frac{1-{\beta }}{{%
\beta }}}\Vert y\Vert _{{\beta }}^{\frac{1}{{\beta }}}.
\end{eqnarray*}%
This proves the inequality (\ref{2.4}).

Assume now that $f$ satisfies (\ref{2.6}). In this case, instead of (\ref%
{2.3}) we have 
\[
|D_{s+}^{\alpha }f(x_{r})|\leq k\left[ \left( a_{0}+a_{1}|x_{r}|\right)
(r-s)^{-\alpha }+\left\| f^{\prime }\right\| _{\infty }\Vert x\Vert _{s,t,{%
\beta }}(r-s)^{\beta -\alpha }\right] . 
\]%
As a consequence,%
\[
\Vert x\Vert _{s,t,{\beta }}\leq k\Vert y\Vert _{{\beta }}\left[
a_{0}+a_{1}\left\| x\right\| _{s,t,\infty }+\left\| f^{\prime }\right\|
_{\infty }\Vert x\Vert _{s,t,{\beta }}(t-s)^{{\beta }}\right] \,.  
\]%
Or%
\[
\Vert x\Vert _{s,t,{\beta }}\leq k\Vert y\Vert _{{\beta }}\left(
a_{0}+a_{1}\left\| x\right\| _{s,t,\infty }\right) \left( 1-k\left\|
f^{\prime }\right\| _{\infty }\Vert y\Vert _{{\beta }}(t-s)^{{\beta }%
}\right) ^{-1}. 
\]%
Therefore,%
\begin{eqnarray*}
|x_{t}| &\leq &|x_{s}|+k\Vert y\Vert _{{\beta }}\left( 1-k\left\| f^{\prime
}\right\| _{\infty }\Vert y\Vert _{{\beta }}(t-s)^{{\beta }}\right) ^{-1}\ 
\\
&&\times \left( a_{0}+a_{1}\left\| x\right\| _{s,t,\infty }\right) {(t-s)}%
^{\beta }.
\end{eqnarray*}%
As before, divide the interval $[0,T]$ into $n=T/\Delta $ subintervals and
set $\Delta =t-s$. Denote 
\begin{eqnarray*}
A &=&k\left\| f^{\prime }\right\| _{\infty }\left\| y\right\| _{{\beta }} \\
B &=&ka_{0}\Vert y\Vert _{{\beta }} \\
C &=&ka_{1}\Vert y\Vert _{{\beta }} \\
D &=&(1-(1-A{\Delta }^{\beta })^{-1}C{\Delta }^{{\beta }})^{-1} \\
F &=&DB(1-A{\Delta }^{\beta })^{-1}\Delta ^{\beta }.
\end{eqnarray*}%
We have 
\begin{eqnarray*}
&&\Vert x\Vert _{s,t,\infty }\left[ 1-k\Vert y\Vert _{{\beta }}(1-A{\Delta }%
^{\beta })^{-1}a_{1}{\Delta }^{\beta }\right] \\
&\leq &|x_{s}|+ka_{0}\Vert y\Vert _{{\beta }}(1-A{\Delta }^{\beta
})^{-1}\Delta ^{\beta }.
\end{eqnarray*}%
This implies 
\[
\sup_{0\leq r\leq t}|x_{r}|\leq (1-(1-A{\Delta }^{\beta })^{-1}C{\Delta }%
^{\beta })^{-1}\left[ \sup_{0\leq r\leq s}|x_{r}|+B(1-A{\Delta }^{\beta
})^{-1}\Delta ^{\beta }\right] \,. 
\]%
Or 
\[
\sup_{0\leq r\leq t}|x_{r}|\leq D\sup_{0\leq r\leq s}|x_{r}|+F. 
\]%
Denote 
\[
Z_{n}=\sup_{0\leq r\leq n{\Delta }}|x_{r}|\,,\quad 
\]%
where $n=\frac{T}{\Delta }$. Then 
\[
Z_{n}\leq DZ_{n-1}+F\leq \cdots \leq D^{n}Z_{0}+\sum_{k=0}^{n-1}D^{k}F. 
\]%
This yields 
\begin{eqnarray*}
\sup_{0\leq t\leq T}|x_{t}| &\leq &(1-(1-A{\Delta }^{\beta })^{-1}C{\Delta }%
^{\beta })^{-T/{\Delta }}|x_{0}| \\
&&+\sum_{k=0}^{n-1}(1-(1-A{\Delta }^{\beta })^{-1}C{\Delta }^{\beta
})^{-k-1}B(1-A{\Delta }^{\beta })^{-1}\Delta ^{\beta }\,.
\end{eqnarray*}%
Then we let ${\Delta }$ satisfy 
\[
A{\Delta }^{\beta }\leq 1/3\,,C{\Delta }^{\beta }\leq 1/3,B{\Delta }^{\beta
}\leq 1/3 
\]%
Namely, we take 
\[
{\Delta }=\left( \frac{1}{3\left( A\vee B\vee C\right) }\right) ^{1/{\beta }%
}\,. 
\]%
Then 
\begin{eqnarray*}
\sup_{0\leq t\leq T}|x_{t}| &\leq &2^{T/{\Delta }}(|x_{0}|+1) \\
&\leq &2^{kT\left[ \Vert f^{\prime }\Vert _{\infty }\vee a_{0}\vee a_{1}%
\right] ^{1/\beta }\Vert y\Vert _{0,T,{\beta }}^{1/{\beta }}}(|x_{0}|+1)\,.
\end{eqnarray*}%
The proof of the theorem is now complete.
\end{proof}

Suppose now that we have two differential equations of the form 
\[
x_{t}=x_{0}+\int_{0}^{t}f(x_{s})dy_{s}, 
\]%
and 
\[
\tilde{x}_{t}=\tilde{x}_{0}+\int_{0}^{t}\tilde{f}(\tilde{x}_{s})\tilde{y}%
_{s}\,, 
\]%
where $y$ and $\widetilde{y}$ are H\"{o}lder continuous functions of order $%
\beta >1/2$, and $f$ and $\widetilde{f}$ are two functions which are
continuously differentiable with H\"{o}lder continuous derivatives of order $%
\lambda >\frac{1}{\beta }-1$. Then, we have the following estimate.

\begin{theorem}
Suppose in addition that $f$ is twice continuously differentiable and $%
f^{\prime \prime }$ is bounded. Then there is a constant $k$ such that%
\begin{eqnarray*}
\sup_{0\leq r\leq T}\left| x_{r}-\tilde{x}_{r}\right| &\leq &2^{kD^{1/{\beta 
}}\left\| y\right\| _{0,T,\beta }^{1/\beta }T} \\
&&\times \Bigg\{|x_{0}-\tilde{x}_{0}|+\Vert y\Vert _{0,T,{\beta }}\left[
\Vert f-\tilde{f}\Vert _{\infty }+\Vert x\Vert _{0,T,{\beta }}\Vert
f^{\prime }-\tilde{f}^{\prime }\Vert _{\infty }\right] \\
&&\quad +\left[ \Vert f\Vert _{\infty }+\Vert \tilde{f}\Vert _{\infty }\Vert
x\Vert _{0,T,\infty }\right] \Vert y-\tilde{y}\Vert _{0,T,{\beta }}\Bigg\}
\end{eqnarray*}
where 
\[
D=\Vert f^{\prime }\Vert _{\infty }\vee \left( \Vert f^{\prime }\Vert
_{\infty }\Vert y\Vert _{0,T,{\beta }}+\Vert f^{\prime \prime }\Vert
_{\infty }(\Vert x\Vert _{0,T,{\beta }}+\Vert \tilde{x}\Vert _{0,T,{\beta }%
})T^{\beta }\right) . 
\]
\end{theorem}

\begin{proof}
Fix $s,t\in \lbrack 0,T]$. Set 
\[
x_{t}-\tilde{x}_{t}-\left( x_{s}-\tilde{x}_{s}\right) =I_{1}+I_{2}+I_{3} 
\]%
where 
\begin{eqnarray*}
I_{1} &=&\int_{s}^{t}\left[ f(x_{r})-f(\tilde{x}_{r})\right] dy_{r} \\
I_{2} &=&\int_{s}^{t}\left[ f(\tilde{x}_{r})-\tilde{f}(\tilde{x}_{r})\right]
dy_{r} \\
I_{3} &=&\int_{s}^{t}\tilde{f}(\widetilde{x}_{r})d\left[ y_{r}-\tilde{y}_{r}%
\right] .
\end{eqnarray*}%
The terms $I_{2}$ and $I_{3}$ can be estimated easily. 
\[
\left| I_{2}\right| \leq k\Vert y\Vert _{{\beta }}\left[ \Vert f-\tilde{f}%
\Vert _{\infty }(t-s)^{\beta }+\Vert f^{\prime }-\tilde{f}^{\prime }\Vert
_{\infty }\Vert \tilde{x}\Vert _{s,t,{\beta }}(t-s)^{2\beta }\right] 
\]%
and 
\[
\left| I_{3}\right| \leq k\Vert y-\tilde{y}\Vert _{{\beta }}\left[ \Vert 
\tilde{f}\Vert _{\infty }(t-s)^{\beta }+\Vert \tilde{f}^{\prime }\Vert
_{\infty }\Vert \tilde{x}\Vert _{s,t,{\beta }}(t-s)^{2\beta }\right] \,, 
\]%
where $\Vert y\Vert _{{\beta }}=\Vert y\Vert _{0,T,{\beta }}$ and $\Vert y-%
\tilde{y}\Vert _{{\beta }}=\Vert y-\tilde{y}\Vert _{0,T,{\beta }}$. The term 
$I_{1}$ is a little more complicated. 
\begin{eqnarray*}
\left| I_{1}\right| &\leq &\int_{s}^{t}|D_{s+}^{\alpha }\left[ f(x_{r})-f(%
\tilde{x}_{r})\right] ||D_{t-}^{1-\alpha }y_{t-}(r)|dr \\
&\leq &k\int_{s}^{t}\Vert y\Vert _{s,t,{\beta }}(t-r)^{\alpha +{\beta }-1}%
\left[ |f(x_{r})-f(\tilde{x}_{r})|(r-s)^{-\alpha }  \right.  \\
&&\left.\quad +\Vert f^{\prime }\Vert
_{\infty }\Vert x-\tilde{x}\Vert _{s,r,{\beta }}(r-s)^{{\beta }-\alpha
}\right. \\
&&\quad \left. +\Vert f^{\prime \prime }\Vert _{\infty }\Vert x-\tilde{x}%
\Vert _{s,r,\infty }\left[ \Vert x\Vert _{s,r,{\beta }}+\Vert \tilde{x}\Vert
_{s,r,{\beta }}\right] (r-s)^{{\beta }-\alpha }\right] dr \\
&\leq &k\Vert y\Vert _{{\beta }}\left\{ \Vert f^{\prime }\Vert _{\infty
}\Vert x-\tilde{x}\Vert _{s,t,\infty }(t-s)^{\beta }+\Vert f^{\prime }\Vert
_{\infty }\Vert x-\tilde{x}\Vert _{s,t,{\beta }}(t-s)^{2{\beta }}\right. \\
&&\quad \left. +\Vert f^{\prime \prime }\Vert _{\infty }\Vert x-\tilde{x}%
\Vert _{s,t,\infty }\left[ \Vert x\Vert _{s,t,{\beta }}+\Vert \tilde{x}\Vert
_{s,t,{\beta }}\right] (t-s)^{2{\beta }}\right\} \,.
\end{eqnarray*}%
Therefore 
\begin{eqnarray*}
\Vert x-\tilde{x}\Vert _{s,t,{\beta }} &\leq &k\Vert y\Vert _{{\beta }%
}\left\{ \Vert f^{\prime }\Vert _{\infty }\Vert x-\tilde{x}\Vert
_{s,t,\infty }+\Vert f^{\prime }\Vert _{\infty }\Vert x-\tilde{x}\Vert _{s,t,%
{\beta }}(t-s)^{{\beta }}\right. \\
&&\quad +\Vert f^{\prime \prime }\Vert _{\infty }\Vert x-\tilde{x}\Vert
_{s,t,\infty }\left[ \Vert x\Vert _{s,t,{\beta }}+\Vert \tilde{x}\Vert _{s,t,%
{\beta }}\right] (t-s)^{{\beta }} \\
&&\left. +\Vert f-\tilde{f}\Vert _{\infty }+\Vert f^{\prime }-\tilde{f}%
^{\prime }\Vert _{\infty }\Vert \tilde{x}\Vert _{s,t,{\beta }}(t-s)^{\beta
}\right\} \\
&&+k\Vert y-\tilde{y}\Vert _{{\beta }}\left[ \Vert \tilde{f}\Vert _{\infty
}+\Vert \tilde{f}^{\prime }\Vert _{\infty }\Vert \tilde{x}\Vert _{s,t,{\beta 
}}(t-s)^{\beta }\right] .
\end{eqnarray*}%
Rearrange it to obtain 
\begin{eqnarray*}
\Vert x-\tilde{x}\Vert _{s,t,{\beta }} &\leq &k(1-k\Vert f^{\prime }\Vert
_{\infty }\Vert y\Vert _{{\beta }}(t-s)^{{\beta }})^{-1}\Bigg\{\Vert y\Vert
_{{\beta }}\Bigg[\Vert f^{\prime }\Vert _{\infty }\Vert x-\tilde{x}\Vert
_{s,t,\infty } \\
&&\quad +\Vert f^{\prime \prime }\Vert _{\infty }\Vert x-\tilde{x}\Vert
_{s,t,\infty }\left[ \Vert x\Vert _{s,t,{\beta }}+\Vert \tilde{x}\Vert _{s,t,%
{\beta }}\right] (t-s)^{{\beta }} \\
&&+\Vert f-\tilde{f}\Vert _{\infty }+\Vert f^{\prime }-\tilde{f}^{\prime
}\Vert _{\infty }\Vert \tilde{x}\Vert _{s,t,{\beta }}(t-s)^{\beta }\Bigg] \\
&&+k\Vert y-\tilde{y}\Vert _{{\beta }}\left[ \Vert \tilde{f}\Vert _{\infty
}+\Vert \tilde{f}^{\prime }\Vert _{\infty }\Vert \tilde{x}\Vert _{s,t,{\beta 
}}(t-s)^{\beta }\right] \Bigg\}.
\end{eqnarray*}%
Set $\Delta =t-s$, and $A=k\Vert f^{\prime }\Vert _{\infty }\Vert y\Vert _{{%
\beta }}$. Then%
\begin{eqnarray*}
\Vert x-\tilde{x}\Vert _{s,t,\infty } &\leq &|x_{s}-\tilde{x}_{s}|+\Vert x-%
\tilde{x}\Vert _{s,t,{\beta }}(t-s)^{\beta } \\
&\leq &|x_{s}-\tilde{x}_{s}|+k(1-A\Delta ^{\beta })^{-1}\Delta ^{\beta }%
\Bigg\{\Vert y\Vert _{{\beta }}\Bigg[\Vert f^{\prime }\Vert _{\infty }\Vert
x-\tilde{x}\Vert _{s,t,\infty } \\
&&+\Vert f^{\prime \prime }\Vert _{\infty }\Vert x-\tilde{x}\Vert
_{s,t,\infty }\left[ \Vert x\Vert _{s,t,{\beta }}+\Vert \tilde{x}\Vert _{s,t,%
{\beta }}\right] \Delta ^{{\beta }} \\
&&+\Vert f-\tilde{f}\Vert _{\infty }+\Vert f^{\prime }-\tilde{f}^{\prime }\Vert _{\infty
}\Vert \tilde{x}\Vert _{s,t,{\beta }}\Delta ^{{\beta }}\Bigg] \\
&&+k\Vert y-\tilde{y}\Vert _{{\beta }}\left[ \Vert \tilde{f}\Vert _{\infty
}+\Vert \tilde{f}^{\prime }\Vert _{\infty }\Vert \tilde{x}\Vert _{s,t,{\beta 
}}\Delta ^{{\beta }}\right] \Bigg\}.
\end{eqnarray*}%
Denote%
\[
B=k\Vert y\Vert _{{\beta }}\left( \Vert f^{\prime }\Vert _{\infty }+\Vert
f^{\prime \prime }\Vert _{\infty }(\Vert x\Vert _{0,T,{\beta }}+\Vert \tilde{%
x}\Vert _{0,T,{\beta }})T^{\beta }\right) \,. 
\]%
Then%
\begin{eqnarray*}
\Vert x-\tilde{x}\Vert _{s,t,\infty } &\leq &\left( 1-(1-A\Delta ^{\beta
})^{-1}\Delta ^{\beta }B\right) ^{-1} \\
&&\times \Bigg\{|x_{s}-\tilde{x}_{s}|+k(1-A\Delta ^{\beta })^{-1}\Delta
^{\beta } \\
&&\times \Bigg[\Vert y\Vert _{{\beta }}\left[ \Vert f-\tilde{f}\Vert
_{\infty }+\Vert f^{\prime }-\tilde{f}^{\prime }\Vert _{\infty }\Vert \tilde{x}\Vert
_{s,t,{\beta }}\Delta ^{{\beta }}\right] \\
&&+\Vert y-\tilde{y}\Vert _{{\beta }}\left[ \Vert \tilde{f}\Vert _{\infty
}+\Vert \tilde{f}^{\prime }\Vert _{\infty }\Vert \tilde{x}\Vert _{s,t,{\beta 
}}\Delta ^{{\beta }}\right] \Bigg]\Bigg\}.
\end{eqnarray*}%
Let ${\Delta }$ satisfy 
\[
A{\Delta }^{\beta }\leq 1/3\,,\quad B{\Delta }^{\beta }\leq 1/3 
\]%
Namely, we take 
\[
{\Delta }=\left( \frac{1}{3\left( A\vee B\right) }\right) ^{1/{\beta }}\,. 
\]%
Then%
\[
\Vert x-\tilde{x}\Vert _{s,t,\infty }\leq 2 \left[  |x_{s}-\tilde{x}%
_{s}|+C\Delta ^{\beta }\right] , 
\]%
where%
\begin{eqnarray*}
C &=&\frac{3}{2}k\Bigg[\Vert y\Vert _{{\beta }}\left[ \Vert f-\tilde{f}\Vert
_{\infty }+\Vert f^{\prime }-\tilde{f}^{\prime }\Vert _{\infty }\Vert \tilde{x}\Vert
_{s,t,{\beta }}\Delta ^{{\beta }}\right] \\
&&+\Vert y-\tilde{y}\Vert _{{\beta }}\left[ \Vert \tilde{f}\Vert _{\infty
}+\Vert \tilde{f}^{\prime }\Vert _{\infty }\Vert \tilde{x}\Vert _{s,t,{\beta 
}}\Delta ^{{\beta }}\right] \Bigg].
\end{eqnarray*}%
Applying the above estimate recursively we obtain 
\[
\sup_{0\leq r\leq T}\left| x_{r}-\tilde{x}_{r}\right| \leq 2^{n}\left[
|x_{0}-\tilde{x}_{0}|+C\Delta^\beta \right] \,, 
\]%
where $T=n\Delta $. Or we have 
\begin{eqnarray*}
\sup_{0\leq r\leq T}\left| x_{r}-\tilde{x}_{r}\right| &\leq &2^{k(\Vert
f^{\prime }\Vert _{\infty }\vee \left( \Vert f^{\prime }\Vert _{\infty
}+\Vert f^{\prime \prime }\Vert _{\infty }(\Vert x\Vert _{0,T,{\beta }%
}+\Vert \tilde{x}\Vert _{0,T,{\beta }})T^{\beta }\right) )^{1/{\beta }%
}\left\| y\right\| _{0,T,\beta }^{1/\beta }T} \\
&&\times \Bigg\{|x_{0}-\tilde{x}_{0}|+\Vert y\Vert _{0,T,{\beta }}\left[
\Vert f-\tilde{f}\Vert _{\infty }+\Vert x\Vert _{0,T,{\beta }}\Vert
f^{\prime }-\tilde{f}\Vert _{\infty }\right] \\
&&\quad +\left[ \Vert f\Vert _{\infty }+\Vert \tilde{f}\Vert _{\infty }\Vert
x\Vert _{0,T,\infty }\right] \Vert y-\tilde{y}\Vert _{0,T,{\beta }}\Bigg\}.
\end{eqnarray*}
\end{proof}

 \setcounter{equation}{0}
\section{Stochastic differential equations driven by a fBm}

Let $B=\{B_{t},t\geq 0\}$ be an $m$-dimensional fractional Brownian motion
(fBm) with Hurst parameter $H>1/2$. That is, $B$ is a Gaussian centered
process with the covariance function $E(B_{t}^{i}B_{s}^{j})=R_{H}(t,s)\delta
_{ij}$, where%
\[
R_{H}(t,s)=\frac{1}{2}\left( t^{2H}+s^{2H}-|t-s|^{2H}\right) . 
\]%
Consider the stochastic differential equation%
\begin{equation}
X_{t}=X_{0}+\int_{0}^{t}\sigma (X_{s})dB_{s}.  \label{3.1}
\end{equation}%
This equation has a unique solution (see \cite{Ly} and \cite{NR}) provided $%
\sigma $ is continuously differentiable, and $\sigma ^{\prime }$ is bounded
and H\"{o}lder continuous of order $\lambda >\frac{1}{H}-1$. The stochastic
integral is interpreted as a path-wise Riemann-Stieltjes integral.

Then, using the estimate (\ref{2.5}) in Theorem \ref{t1} we obtain the
following estimate for the solution of Equation (\ref{3.1}), if we choose $%
\beta \in \left( \frac{1}{2},H\right) $. Notice that $\frac{1}{\beta }<2.$%
\begin{equation}   \label{eq1}
\sup_{0\leq t\leq T}|X_{t}|\leq 2^{kT\left( \Vert \sigma ^{\prime }\Vert _{\infty
}\vee |\sigma (0)|\right) \Vert B\Vert _{0,T,{\beta }}^{1/{\beta }}}(|X_{0}|+1)\,. 
\end{equation}

If $\sigma $ is bounded we can make use of the estimate (\ref{2.4}) and we
obtain%
\begin{equation}
\sup_{0\leq t\leq T}|X_{t}|\leq |X_{0}|+kT\Vert \sigma \Vert _{\infty }\Vert
\sigma ^{\prime }\Vert _{\infty }^{\frac{1-{\beta }}{{\beta }}}\Vert B\Vert
_{0,T,{\beta }}^{\frac{1}{{\beta }}}.   \label{eq2}
\end{equation}
 
These estimates improve those obtained by Nualart and R\u{a}\c{s}canu in \cite{NR}  
based on a suitable version of Gronwall's lemma.
 The estimates (\ref{eq1}) and (\ref{eq2}) allow us to establish the following integrability
 properties for the solution of Equation (\ref{3.1}).
 
 \begin{theorem}
 Consider the stochastic differential equation (\ref{3.1}). 
If $\sigma'$ is bounded and H\"older continuous 
of order $\lambda >\frac1 H-1$,  then 
\begin{equation}
E\left( \sup_{0\leq t\leq T}|X_{t}|^{p}\right) <\infty 
\end{equation}
for all $p\geq 2$.   If furthermore $\sigma$ is bounded, then 
\begin{equation}
E\left( \exp \lambda \left( \sup_{0\leq t\leq T}|X_{t}|^{\gamma }\right)
\right) <\infty 
\end{equation}
for any $\lambda >0$ and $\gamma <2\beta $.
 \end{theorem}
 
If we  apply these results to the linear equation satisfied by the
derivative in the sense of Malliavin calculus of $X_{t}$  then we get that $%
X_{t}$ belongs to the Sobolev space $\mathbb{D}^{1,p}$ for all $p\geq 2$.  This
implies that  if the coefficient $\sigma $ is infinitely differentiable with
bounded derivatives of all orders, then, $X_{t}$ belongs to $\mathbb{D}%
^{\infty }$.   This allows us to deduce the regularity of the density of the
random vector $X_{t}$ at a fixed time $t>0$ assuming the following
nondegeneracy condition:

(H) The vector space spanned by $\left\{ \left( \sigma ^{ij}(X_{0})\right)
_{1\leq i\leq d},1\leq j\leq m\right\} $ is $\mathbb{R}^{m}$.

That is, we have the following result.

\begin{theorem} Consider the stochastic differential equation (\ref{3.1}). 
 Suppose that $\sigma $ is infinitely differentiable with bounded derivatives
of all orders, and  the assumption (H) holds. Then, for any $t>0$ the probability law of  $X_t$ 
has an  $C^\infty$ density. 
\end{theorem}
 
In \cite{NS} Nualart and Saussereau have proved that the random variable $X_t$
belongs locally to the space $\mathbb{D}^{\infty }$, and, as  a consequence, they have derived the absolute continuity of the law of $X_t$ under the assumption (H).

\end{document}